\documentclass[11pt]{article}
\usepackage{amsthm}
\usepackage{amsmath}  
\usepackage{amssymb}  
\usepackage{geometry}  
\usepackage{hyperref}  
\usepackage{authblk}  
\usepackage{setspace}  
\usepackage{pxfonts}
\usepackage{yfonts}
\usepackage{dsfont}
\usepackage{graphicx}
\usepackage{relsize}
\usepackage{color}
\usepackage{marvosym}
\usepackage{graphicx}
\usepackage{relsize}
\usepackage[latin1]{inputenc}
\usepackage{enumerate}

\def\bel{\begin{equation}\label}
\def\eeq{\end{equation}}
\def\ds{\displaystyle}

\def\mt{\longrightarrow}
\def\v{\vskip 1em}

\def\R{\mathbb R}
\def\Z{{\bf Z}}
\def\W{{\bf W}}
\def\C{\mathfrak{B}}

\def\N{{\bf N}}

\def\S{{\bf S}}

\def\Q{{\bf Q}}

\def\H{{\bf H}}
\def\L{{\bf L}}
\def\U{{\bf U}}

\def\T{{\bf T}}

\def\p{{\partial}}
\def\a{{\bf a}}
\def\b{{\bf b}}
\def\e{{\bf e}}

\def\supp{{\bf supp}}
\def\vol{{\bf vol}}

\def\I{{\bf I}}

\def\alpha{\alphaup}
\def\beta{\betaup}
\def\gamma{\gammaup}
\def\delta{\deltaup}
\def\theta{\thetaup}
\def\xi{{\xiup}}
\def\eta{{\etaup}}
\def\tau{{\tauup}}
\def\rho{{\rhoup}}
\def\phi{{\phiup}}
\def\psi{{\psiup}}
\def\lambda{{\lambdaup}}
\def\omega{\omegaup}
\def\varphi{{\varphiup}}
\def\gamma{{\gammaup}}
\def\c{{\bf c}}

\newtheorem{remark}{Remark}[section]

\begin{document}

\[\begin{array}{cc}\ds\hbox{\LARGE{\bf $\H^1\mt\L^q$-boundedness of fractional integral operators}}
\\\\ \ds
\hbox{\LARGE{\bf  having a flag kernel}}\end{array}\]

\[\hbox{Jiashu Zhang}\qquad\hbox{and}\qquad\hbox{Zipeng Wang} \]
\begin{abstract}
We study a family of fractional integral operators
\[\I_{\alpha\beta}^\rho f(x,y)~=~\iint_{\R^n\times\R^m} f(x-u,y-v) \left({1\over |u|}\right)^{n-\alpha} \left[{1\over |u|^\rho+|v|}\right]^{m-\beta}dudv\]
for $0<\alpha<n,0<\beta<m$ and $\rho\ge1$.  First, we show  $\I_{\alpha\beta}^\rho\colon \L^p(\R^{n+m})\mt\L^q(\R^{n+m}),1<p<q<\infty$ if and only if ${\alpha\over n}\ge{\beta\over m}$ and ${\alpha+\rho\beta\over n+\rho m}={1\over p}-{1\over q}$. Second, we prove that $\I_{\alpha\beta}^\rho$ is bounded from the classical, atom decomposable $\H^1$-Hardy space to $\L^q(\R^{n+m})$ if and only if ${\alpha\over n}>{\beta\over m}$ and ${\alpha+\rho\beta\over n+\rho m}=1-{1\over q}$.
\end{abstract}

\section{Introduction}
\setcounter{equation}{0}
Let $0<\alpha<n$. A fractional integral operator $\I_\alpha$ is initially defined  as
\bel{I_a}
\I_\alpha f(x)~=~\int_{\R^n}f(u) \left[{1\over|x-u|}\right]^{n-\alpha} du.
\eeq
In 1928, Hardy and Littlewood \cite{Hardy-Littlewood} have obtained an regularity theorem for $\T_\a$ when $\N=1$. Ten years later, Sobolev \cite{Sobolev} made extensions on every higher dimensional space.

$\diamond$ {\small Throughout, $\C>0$ is regarded as a generic constant depending on its sub-indices.}
\v
{\bf Hardy-Littlewood-Sobolev theorem}~~{\it Let $\I_\alpha$ defined in (\ref{I_a}) for $0<\alpha<n$. We have
\bel{HLS Ineq}
\begin{array}{cc}\ds
\left\| \I_\alpha f\right\|_{\L^q(\R^n)}~\leq~\C_{p~q}~\left\| f\right\|_{\L^p(\R^n)},\qquad 1<p<q<\infty
\\\\ \ds
\hbox{if and only if}\qquad
{\alpha\over n}~=~{1\over p}-{1\over q}.
\end{array}
\eeq
}\\
This classical result has been re-investigated by Krantz \cite{Krantz} on Hardy spaces.
\v
{\bf Krantz theorem}~~{\it Let $\I_\alpha$ defined in (\ref{I_a}) for $0<\alpha<n$. We have
\bel{K Ineq}
\begin{array}{cc}\ds
\left\| \I_\alpha f\right\|_{\H^q(\R^n)}~\leq~\C_{p~q}~\left\| f\right\|_{\H^p(\R^n)},\qquad 0<p<q<\infty
\\\\ \ds
\hbox{if and only if}\qquad
{\alpha\over n}~=~{1\over p}-{1\over q}.
\end{array}
\eeq
}

\begin{remark}  $\H^p(\R^{n}), 0<p\leq1$  is the classical $\H^p$-Hardy space  investigated by Fefferman and Stein \cite{Fefferman-Stein}. Moreover, it has a characterization of atomic decomposition established by Coifman \cite{Coifman}.
\end{remark}

Recently, {\bf Krantz theorem} is extended to the multi-parameter setting by Tang \cite{Tang} whereas $\I_\alpha$ in (\ref{I_a}) is replaced with strong fractional integral operators whose kernels have singularity on every coordinate subspace. To better illustrate the difference between this new result and ours, we focus on its bi-parameter version.
Let $0<\alpha<n$, $0<\beta<m$. Define
\bel{I ab}
\I_{\alpha\beta} f(x,y)~=~\iint_{\R^n\times\R^m} f(u,v) \left[{1\over|x-u|}\right]^{n-\alpha} \left[{1\over|y-v|}\right]^{m-\beta}dudv.
\eeq

{\bf Tang theorem}~~{Let $\I_{\alpha\beta}$ defined in (\ref{I ab}) for $0<\alpha<n$, $0<\beta<m$. We have
\bel{T Ineq}
\begin{array}{cc}\ds
\left\| \I_{\alpha\beta} f\right\|_{\H^q\times\H^q(\R^n\times\R^m)}~\leq~\C_{p~q}~\left\| f\right\|_{\H^p\times\H^p(\R^n\times\R^m)},\qquad 0<p<q<\infty
\\\\ \ds
\hbox{if and only if}\qquad
{\alpha\over n}~=~{\beta\over m}~=~{1\over p}-{1\over q}.
\end{array}
\eeq
}

\begin{remark} $\H^p\times \H^p(\R^n\times\R^m), 0<p\leq1$ is  the product Hardy space introduced by   Gundy and Stein \cite{Gundy-Stein}. Furthermore, it  cannot be characterized in terms of  "rectangle atoms". See the counter-example of Carleson \cite{Carleson}. 
\end{remark}

In this paper, we study a family of fractional integral operators  
whose kernels  satisfying non-isotropic dilations have singularity on a coordinate subspace, commonly known as flag kernels.  An initial motivation to assert these operators  comes from certain  sub-elliptic boundary value problems.  The  solution turns out to be a composition of two singular integral operators. One of them is elliptic. The other is parabolic associated with an non-isotropic dilation. Singular integrals of this type have been systematically studied. For instance, see the paper by Phong and Stein \cite{Phong-Stein},
Muller, Ricci and Stein \cite{Muller-Ricci-Stein}, Nagel, Ricci and Stein \cite{Nagel-Ricci-Stein}, Nagel, Ricci, Stein and Wainger \cite{Nagel-Ricci-Stein-Wainger}, Han-et-al \cite{Han-Lin-Lu-Ruan-Sawyer}, Han, Lu and Sawyer \cite{Han-Lu-Sawyer} and Han, Lee, Li and Wick \cite{Han-Lee-Li-Wick'}-\cite{Han-Lee-Li-Wick}. This direction remains largely open for fractional integrals.

Let $0<\alpha<n$, $0<\beta<m$ and $\rho\ge1$. We define
\bel{I abrho}
\I_{\alpha\beta}^\rho f(x,y)~=~\iint_{\R^n\times\R^m} f(u,v) \Omega^{\alpha\beta}_\rho(x-u,y-v)dudv
\eeq
where $\Omega^{\alpha\beta}_\rho$ is a distribution on $\R^{n+m}$ agree with
\bel{Omega}
\Omega^{\alpha\beta}_\rho(x,y)~=~\left({1\over |x|}\right)^{n-\alpha} \left[{1\over |x|^\rho+|y|}\right]^{m-\beta},\qquad\hbox{\small{$x\neq0$}}.
\eeq
A concrete example is given by Muller, Ricci and Stein \cite{Muller-Ricci-Stein}. Let  $x=(z,w)\in\R^d\times\R^d$. Consider $\mathcal{L}^{-\a}T^{-\b}$ for $0<\a<d, 0<\b<1$:  $T=\partial_t$ and $\mathcal{L}$ is the sub-Laplacian:
$
\mathcal{L}=-\sum_{i=1}^d \Z_i^2+\W^2_i,~ \Z_i=\p_{z_i}+2 w_i \p_t,~ \W_i=\p_{w_i}-2z_i \p_t$.
The inverse of $\mathcal{L}^\a$  is given as the Riesz potential:
\[
\mathcal{L}^{-\a}~=~{1\over \Gamma(\a)}\int_0^\infty s^{\a-1} e^{-s\mathcal{L}} ds
\]
where $\Gamma$ denotes Gamma function. 
The kernel of $\mathcal{L}^{-\a}T^{-\b}$ is a distribution in $\R^{2d+1}$ agree with a function similar to
$ \Gamma\Big({1-\b\over 2}\Big) \left({1\over |x|}\right)^{2d-2\a} \left[{1\over |x|^2+|y|}\right]^{1-\b}$ for $x\neq0$.\footnote{We say $A$ similar to $B$ if $\c^{-1} B\leq A\leq \c B$ for some $\c>0$. } In compare to (\ref{Omega}), we find $n=2d, m=1$, $\alpha=2\a, \beta=\b$ and $\rho=2$.

For $1<p<q<\infty$, we have  $\mathcal{L}^{-\a}T^{-\b}\colon \L^p(\R^{2d+1})\mt \L^q(\R^{2d+1})$ if and only if $\a \ge d\b$ and ${\a+\b\over d+1}={1\over p}-{1\over q}$. This  is proved in section 6 of \cite{Muller-Ricci-Stein} by using complex interpolation. One of the two end-point estimates relies on the $\L^p$-theorem developed there. First, we  show that every convolution operator with a kernel similar to (\ref{Omega}) satisfies the desired $\L^p\mt\L^q$-regularity.

\v
{\bf Theorem One}~~{\it Let $\I_{\alpha\beta}^\rho$ defined in (\ref{I abrho})-(\ref{Omega}) for $0<\alpha<n$, $0<\beta<m$ and $\rho\ge1$. We have
\bel{Result One}
\left\| \I_{\alpha\beta}^\rho f\right\|_{\L^q(\R^{n+m})}~\leq~\C_{\alpha~\beta~\rho~p~q} ~\left\| f\right\|_{\L^p(\R^{n+m})},\qquad 1<p<q<\infty
\eeq
if and only if
\bel{Formula One}
{\alpha\over n}\ge{\beta\over m},\qquad {\alpha+\rho\beta\over n+\rho m}~=~{1\over p}-{1\over q}.
\eeq}\\
Our main result is to give a characterization for the $\H^1\mt\L^q$-boundedness of $\I_{\alpha\beta}^\rho$.
\v
{\bf Theorem Two}~~{\it Let $\I_{\alpha\beta}^\rho$ defined in (\ref{I abrho})-(\ref{Omega}) for $0<\alpha<n$, $0<\beta<m$ and $\rho\ge1$. We have
\bel{Result Two}
\left\| \I_{\alpha\beta}^\rho f\right\|_{\L^q(\R^{n+m})}~\leq~\C_{\alpha~\beta~\rho~q} ~\left\| f\right\|_{\H^1(\R^{n+m})},\qquad 1<q<\infty
\eeq
if and only if
\bel{Formula Two}
{\alpha\over n}>{\beta\over m},\qquad {\alpha+\rho\beta\over n+\rho m}~=~1-{1\over q}.
\eeq}

Note that $\H^1(\R^{n+m})$ is the same $\H^1$-Hardy space introduced by Fefferman and Stein \cite{Fefferman-Stein} which has an atomic decomposition due to Coifman \cite{Coifman}. On the other hand, the product Hardy space $\H^1\times\H^1(\R^n\times\R^m)$ defined by Gundy and Stein \cite{Gundy-Stein} is a subspace of $\H^1(\R^{n+m})$.

The fractional integral operator $\I_{\alpha\beta}^\rho$ whose  kernel $\Omega^{\alpha\beta}_\rho$ carries certain multi-parameter structure  as defined in (\ref{Omega}) is still bounded from the classical, atom decomposable $\H^1$-Hardy space to $\L^q(\R^{n+m})$.
\v

The remaining paper is organized as follows. In the next section, we show $\I_{\alpha\beta}^\rho\colon\L^p(\R^{n+m})\mt\L^q(\R^{n+m})$  implying ${\alpha\over n}\ge{\beta\over m}$ and ${\alpha+\rho\beta\over n+\rho m}={1\over p}-{1\over q}$ for $1\leq p<q<\infty$. Moreover, we give a counter example for $ \I_{\alpha\beta}^\rho\colon\H^1(\R^{n+m})\mt\L^q(\R^{n+m})$ when ${\alpha\over n}={\beta\over m}$ and ${\alpha+\rho\beta\over n+\rho m}=1-{1\over q}$. In section 3, we show that (\ref{Formula One}) implies (\ref{Result One}).
Section 4 is devoted to the proof of (\ref{Formula Two}) $\Longrightarrow$ (\ref{Result Two}). By the characterization of atomic decomposition for $\H^1(\R^{n+m})$, it is suffice to show $\left\|\I_{\alpha\beta}^\rho a\right\|_{\L^q(\R^{n+m})}\leq\C_{\alpha~\beta~\rho~q}$ of which $a$ is an $\H^1$-atom:
\bel{atom}
\supp a\subset \Q,\qquad |a(x,y)|\leq {1\over \vol\{\Q\}},\qquad \iint_{\Q} a(x,y)dxdy~=~0
\eeq
where $\Q\subset\R^{n+m}$ is some cube parallel to the coordinates.

\section{Proof of necessary conditions}
\setcounter{equation}{0}
Let $\I_{\alpha\beta}^\rho$ defined in (\ref{I abrho})-(\ref{Omega}) for $0<\alpha<n$, $0<\beta<m$ and $\rho\ge1$.  We have
\bel{I rewrite}
\begin{array}{lr}\ds
\I_{\alpha\beta}^\rho f(x,y)~=~\iint_{\R^n\times\R^m} f(u,v) \left({1\over |x-u|}\right)^{n-\alpha} \left[{1\over |x-u|^\rho+|y-v|}\right]^{m-\beta}dudv.
\end{array}
\eeq
By changing dilations $(x,y)\mt (\delta x, \delta^\rho \lambda y)$ and $(u,v)\mt (\delta u, \delta^\rho\lambda v)$ for $\delta>0, \lambda>1$,  we find
\bel{Dila I}
\begin{array}{lr}\ds
\left\{ \iint_{\R^n\times\R^m} 
\left\{\iint_{\R^n\times\R^m} f\left[\delta^{-1} u,\delta^{-\rho}\lambda^{-1}v\right]\left({1\over |x-u|}\right)^{n-\alpha}\left[{1\over |x-u|^\rho+|y-v|}\right]^{m-\beta}
 du dv \right\}^q dxdy\right\}^{1\over q}
 \\\\ \ds
=~\left\{ \iint_{\R^n\times\R^m} 
\left\{\iint_{\R^n\times\R^m} f( u,v) \left({1\over \delta|x-u|}\right)^{n-\alpha}\left[{1\over \delta^\rho|x-u|^\rho+\delta^\rho \lambda|y-v|}\right]^{m-\beta}
\delta^{n+\rho m}\lambda^m du dv \right\}^q \delta^{n+\rho m}\lambda^m dxdy\right\}^{1\over q} 
 \\\\ \ds
=~\delta^{\alpha+\rho\beta} \delta^{n+\rho m\over q} \lambda^{m\over q}
\left\{ \iint_{\R^n\times\R^m} 
\left\{\iint_{\R^n\times\R^m} f( u,v) \left({1\over |x-u|}\right)^{n-\alpha}\left[{1\over |x-u|^\rho+\lambda|y-v|}\right]^{m-\beta}
\lambda^m du dv \right\}^q dxdy\right\}^{1\over q}
 \\\\ \ds
\ge~\delta^{\alpha+\rho\beta} \delta^{n+\rho m\over q} \lambda^{m\over q}
\\ \ds~~~
\left\{ \iint_{\R^n\times\R^m} 
\left\{\iint_{\R^n\times\R^m} f( u,v) \left({1\over |x-u|}\right)^{n-\alpha}\left[{1\over |x-u|^\rho+|y-v|}\right]^{m-\beta}
\lambda^{\beta-m}\lambda^m du dv \right\}^q dxdy\right\}^{1\over q}\quad \hbox{\small{($\lambda>1$)}}
\\\\ \ds
=~\delta^{\alpha+\rho\beta} \delta^{n+\rho m\over q} \lambda^\beta \lambda^{m\over q}
\left\{ \iint_{\R^n\times\R^m} 
\left\{\iint_{\R^n\times\R^m} f( u,v) \left({1\over |x-u|}\right)^{n-\alpha}\left[{1\over |x-u|^\rho+|y-v|}\right]^{m-\beta}
 du dv \right\}^q dxdy\right\}^{1\over q}.
\end{array}
\eeq
Consider
\bel{Norm Ineq}
\left\| \I_{\alpha\beta}^\rho f\right\|_{\L^q(\R^{n+m})}~\lesssim~\left\| f\right\|_{\L^p(\R^{n+m})},\qquad 1\leq p<q<\infty
\eeq
which implies that the last line of (\ref{Dila I}) is bounded by a constant multiple of
\bel{I L^p}
\begin{array}{lr}\ds
\left\{ \iint_{\R^n\times\R^m} 
\left\{\iint_{\R^n\times\R^m} f\left[\delta^{-1} u,\delta^{-\rho}\lambda^{-1}v\right]\left({1\over |x-u|}\right)^{n-\alpha}\left[{1\over |x-u|^\rho+|y-v|}\right]^{m-\beta}
 du dv \right\}^q dxdy\right\}^{1\over q}
\\\\ \ds
~\lesssim~\left\{ \iint_{\R^n\times\R^m} \Big[ f\left(\delta^{-1} x,\delta^{-\rho}\lambda^{-1} y\right) \Big]^p dx dy\right\}^{1\over p}
\\\\ \ds
~=~\delta^{n+\rho m\over p} \lambda^{m\over p} \left\| f\right\|_{\L^p(\R^{n+m})}.
\end{array}
\eeq
This must be true for every $\delta>0$ and $\lambda>1$. We necessarily have
\bel{homogeneity}
{\alpha+\rho\beta\over n+\rho m}~=~{1\over p}-{1\over q}
\eeq
and
\bel{beta<}
\beta~\leq~{m\over p}-{m\over q}.
\eeq
By putting together (\ref{homogeneity}) and (\ref{beta<}), we find
\bel{am>bn}
\alpha m~\ge~\beta n.
\eeq

\subsection{A counter example for $\H^1\mt\L^q$-estimate}
Now, we give a counterexample for $\I_{\alpha\beta}^\rho\colon\H^1(\R^{n+m})\mt\L^q(\R^{n+m})$ when  $\alpha m=n\beta$. 

Let $\Q_o=\{(x,y)\in\R^n\times\R^m\colon |x_i|\leq1, i=1,2,\ldots,n;~~|y_i|\leq1, i=1,2,\ldots,m\}$. Consider $a(x,y)=\operatorname{sgn}(x_1)\chi_{\Q_o}(x,y)$ which is an $\H^1$-atom in $\R^{n+m}$. Define $\U\subset\R^n$ by $\U=\big\{x\in\R^n\colon 2\leq x_i\leq 4, i=1,2,\ldots,n\big\}$. We aim to show 
\bel{Int infinity}
\iint_{\U\times\R^m} \left|\I_{\alpha\beta}^\rho a(x,y)\right|^qdxdy~=~\infty.
\eeq
Denote $\e_1$ to be the unit vector in $x_1$-axis and $\Q^+_o=\Q_o\cap\{x_1>0\}$. For $(x,y)\in\U\times\R^m$, we find
\bel{I > est}
\begin{array}{lr}\ds
\I_{\alpha\beta}^\rho a(x,y)~=~\iint_{\Q_o} a(u,v)\Omega^{\alpha\beta}_\rho (x-u,y-v)dudv
\\\\ \ds~~~~~~~~~~~~~~
~=~\iint_{\Q^+_o} \left({1\over |x-u|}\right)^{n-\alpha} \left[{1\over |x-u|^\rho+|y-v|}\right]^{m-\beta} -  \left({1\over |x-u+\e_1|}\right)^{n-\alpha} \left[{1\over |x-u+\e_1|^\rho+|y-v|}\right]^{m-\beta}dudv
\\\\ \ds~~~~~~~~~~~~~~
~\gtrsim~\iint_{\Q^+_o} \left[\frac{1}{1+|y-v|}\right]^{m-\beta}dudv
\\\\ \ds~~~~~~~~~~~~~~
~\gtrsim~ \left(\frac{1}{1+|y|}\right)^{m-\beta}.
\end{array}
\eeq
Note that ${\alpha\over n}={\beta\over m}$ and ${\alpha+\rho\beta\over n+\rho m}=1-{1\over q}$ together imply 
$\beta=m-{m\over q}$. We have
\bel{Int > infty}
\begin{array}{lr}\ds
\iint_{\U\times\R^m}\left|\I_{\alpha\beta}^\rho a(x,y)\right|^q dxdy ~\gtrsim~\iint_{\U\times\R^m} \left(\frac{1}{1+|y|}\right)^{q(m-\beta)} dxdy\qquad\hbox{\small{by (\ref{I > est})}}
\\\\ \ds~~~~~~~~~~~~~~~~~~~~~~~~~~~~~~~~~~~~~~~~~
~\approx~ \int_{\R^m}\left(\frac{1}{1+|y|}\right)^m dy ~=~\infty.
\end{array}
\eeq

\section{Proof of Theorem One}
\setcounter{equation}{0}
Recall (\ref{Formula One}). We have ${\alpha\over n}\ge{\beta\over m}$ and $ {\alpha+\rho\beta\over n+\rho m}={1\over p}-{1\over q}$ for $1<p<q<\infty$.
Define $0<\a\leq\alpha<n$ and $0<\beta\leq\b<m$ implicitly by requiring 
\bel{ab}
{\a\over n}~=~{\b\over m},\qquad \a+\rho \b~=~\alpha+\rho\beta.
\eeq
By solving the two equations in (\ref{ab}), we find 
\bel{ab formula}
\a={\alpha+\rho\beta\over 1+\rho({m\over n})},\qquad \b={\alpha+\rho\beta\over {n\over m}+\rho}.
\eeq
Furthermore, (\ref{ab}) together with the homogeneity condition $ {\alpha+\rho\beta\over n+\rho m}={1\over p}-{1\over q}$ imply
\bel{Formula ab new}
{\a\over n}~=~{\b\over m}~=~{\a+\rho\b\over n+\rho m}~=~{\alpha+\rho\beta\over n+\rho m}~=~{1\over p}-{1\over q}.
\eeq

Let $\I_{\alpha\beta}^\rho f=f\ast\Omega^{\alpha\beta}_\rho$ defined in (\ref{I abrho})-(\ref{Omega}). Observe that
\bel{Omega dominate}
\begin{array}{lr}\ds
\Omega^{\alpha\beta}_\rho(x,y)~=~\left({1\over |x|}\right)^{n-\alpha} \left[{1\over |x|^\rho+|y|}\right]^{m-\beta}
\\\\ \ds~~~~~~~~~~~~~~~
~=~\left({1\over |x|}\right)^{n-\alpha} \left[{1\over |x|^\rho+|y|}\right]^{m-\b} \left[{1\over |x|^\rho+|y|}\right]^{\b-\beta}
\\\\ \ds~~~~~~~~~~~~~~~
~\leq~\left({1\over |x|}\right)^{n-\alpha} \left[{1\over |x|^\rho+|y|}\right]^{m-\b} \left({1\over |x|}\right)^{\rho\b-\rho\beta}
\\\\ \ds~~~~~~~~~~~~~~~
~=~\left({1\over |x|}\right)^{n-\a} \left[{1\over |x|^\rho+|y|}\right]^{m-\b}
\\\\ \ds~~~~~~~~~~~~~~~
~\leq~\left({1\over |x|}\right)^{n-\a} \left({1\over |y|}\right)^{m-\b}\qquad\hbox{\small{$x\neq0$}}.
\end{array}
\eeq
Recall {\bf Hardy-Littlewood-Sobolev theorem} stated in the beginning of this paper. 
We have
\bel{Iteration}
\begin{array}{lr}\ds
\left\| \I_{\alpha\beta}^\rho f\right\|_{\L^q(\R^{n+m})}~=~
\left\{\iint_{\R^n\times\R^m} \left\{\iint_{\R^n\times\R^m} f(u,v) \Omega^{\alpha\beta}_\rho(x-u,y-v)dudv\right\}^q dxdy\right\}^{1\over q}
\\\\ \ds~~~~~~~~~~~~~~~~~~~~~~
~\leq~\left\{\iint_{\R^n\times\R^m} \left\{\iint_{\R^n\times\R^m} f(u,v)\left({1\over |x-u|}\right)^{n-\a} \left({1\over |y-v|}\right)^{m-\b}dudv\right\}^q dxdy\right\}^{1\over q}\qquad\hbox{\small{ by (\ref{Omega dominate})}}
\\\\ \ds~~~~~~~~~~~~~~~~~~~~~~
~\leq~\left\{\int_{\R^m} \left\{\int_{\R^n}\left\{\int_{\R^m} f(x,v) \left({1\over |y-v|}\right)^{m-\b}dv\right\}^p dx\right\}^{q\over p} dy\right\}^{1\over q}\qquad\hbox{\small{ by (\ref{Formula ab new}) and (\ref{HLS Ineq})}}
\\\\ \ds~~~~~~~~~~~~~~~~~~~~~~
~\leq~\left\{\int_{\R^n} \left\{\int_{\R^m}\left\{\int_{\R^m} f(x,v) \left({1\over |y-v|}\right)^{m-\b}dv\right\}^q dy\right\}^{p\over q} dx\right\}^{1\over p}~~~~\hbox{\small{ by Minkowski intergal inequality}}
\\\\ \ds~~~~~~~~~~~~~~~~~~~~~~
~\leq~\left\| f\right\|_{\L^p(\R^{n+m})}\qquad \hbox{\small{ by (\ref{Formula ab new}) and (\ref{HLS Ineq})}}.
\end{array}
\eeq

\section{Proof of Theorem Two}
\setcounter{equation}{0}
 Let $0<\alpha<n$, $0<\beta<m$ and $\rho\ge1$. From (\ref{Omega}), we have 
\[\Omega^{\alpha\beta}_\rho(x,y)~=~\left({1\over |x|}\right)^{n-\alpha} \left[{1\over |x|^\rho+|y|}\right]^{m-\beta},\qquad\hbox{\small{$x\neq0$}}.\]
A direct computation gives
\bel{deri est}
\begin{array}{cc}\ds
\left|\nabla_x\Omega^{\alpha\beta}_\rho(x,y)\right|~\leq~\C_{\alpha~\beta~\rho}~ \Omega^{\alpha\beta}_\rho(x,y)~\max\left\{{1\over|x|},{|x|^{\rho-1}\over |x|^{\rho}+|y|}\right\},
\\\\ \ds
\left|\nabla_y\Omega^{\alpha\beta}_\rho(x,y)\right|~\leq~\C_{\alpha~\beta~\rho}~ \Omega^{\alpha\beta}_\rho(x,y)~{1\over |x|^{\rho}+|y|}.
\end{array}
\eeq
From (\ref{deri est}), we conclude
\bel{deri Est}
\left|\nabla \Omega^{\alpha\beta}_\rho(x,y)\right|~\leq~\C_{\alpha~\beta~\rho}~ \Omega^{\alpha\beta}_\rho(x,y)\max\left\{{1\over|x|},{1\over |x|^{\rho}+|y|}\right\}.
\eeq
Let $\I_{\alpha\beta}^\rho f=f\ast\Omega^{\alpha\beta}_\rho$ defined in (\ref{I abrho})-(\ref{Omega}).
We aim to prove		
\bel{Result Two a}
\left\| \I_{\alpha\beta}^\rho a\right\|_{\L^q(\R^{n+m})}~\leq~ \C_{\alpha~\beta~\rho~q}\qquad\hbox{for}\qquad {\alpha\over n}>{\beta\over m},\quad\frac{\alpha+\rho\beta}{n+\rho m}~=~1-{1\over q}	
\eeq	
where $a$ is an $\H^1$-atom in $\R^{n+m}$ satisfying (\ref{atom}).

Denote $\Q\subset\R^{n+m}$ to be a cube centered on the origin with a side length equal to $2^L$.
Without lose of the generality, we assume $a$ supported in the cube ${1\over 2}\Q$: co-center with $\Q$ having a side length $2^{L-1}$.
For $k,\ell\ge0$, we define 
		\bel{Q_kl}
\Q=\Q_{00},\qquad	\Q_{k\ell}=\Bigg\{(x,y)\in\R^n\times\R^m\colon2^{L+k-1}\leq|x|< 2^{L+k},~~2^{L+\ell-1}\leq|y|< 2^{L+\ell}\Bigg\}.
			\eeq
	\begin{figure}[h]
\centering
\includegraphics[scale=0.49]{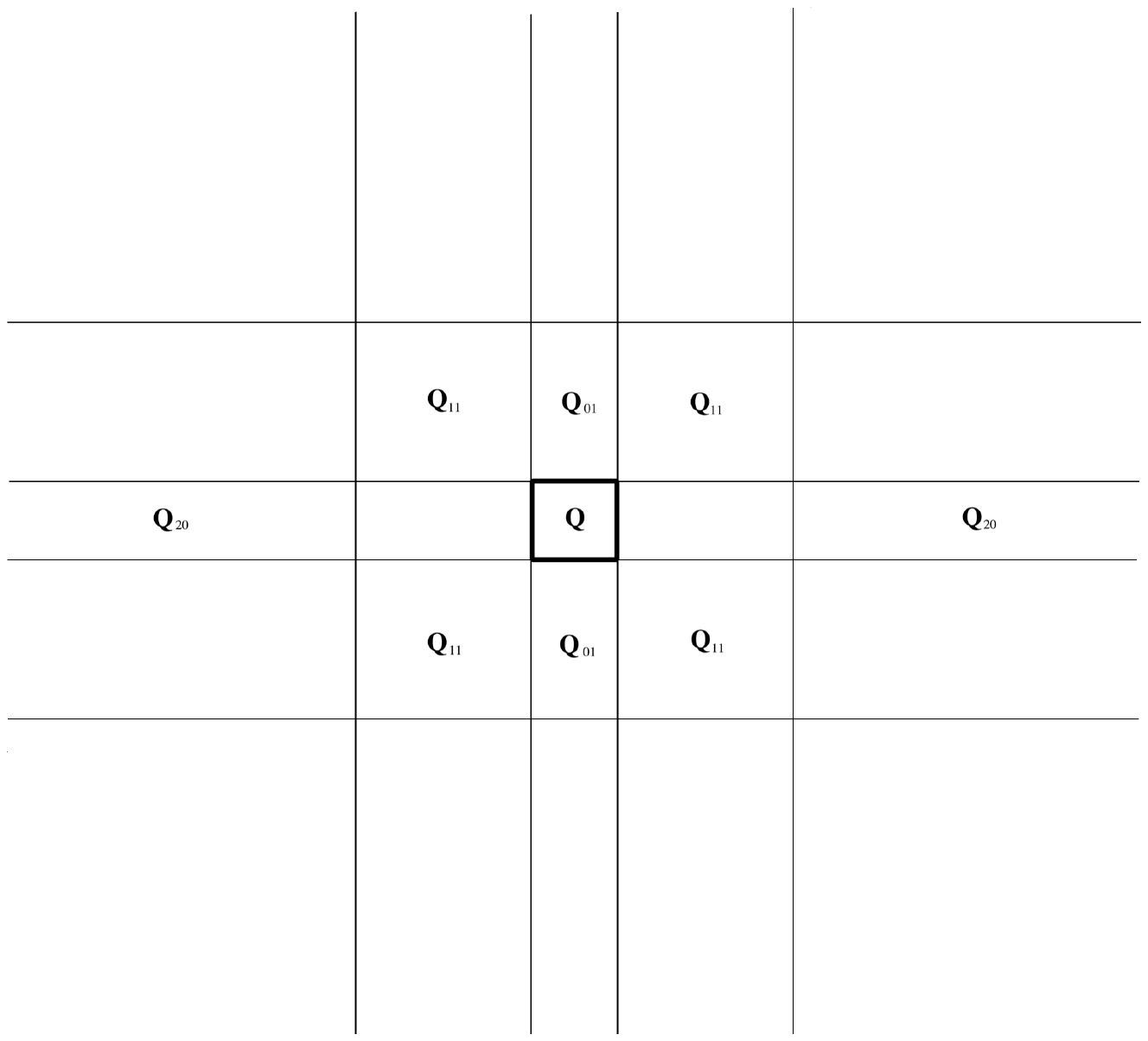}
\end{figure}

We obtain (\ref{Result Two a}) by estimating
\bel{Result Two a}
\iint_{\Q_{k\ell}} \left|\I_{\alpha\beta}^\rho a(x,y)\right|^q dxdy,\qquad {\alpha\over n}>{\beta\over m},\quad \frac{\alpha+\rho\beta}{n+\rho m}~=~1-{1\over q}
\eeq				
$w.r.t$ {\bf Case 1}: $k=\ell=0$, {\bf Case 2}: $k>0,\ell>0$, {\bf Case 3}: $k>0,\ell=0$ and {\bf Case 4}: $k=0,\ell>0$.

Before moving forward, we note that ${\alpha\over n}>{\beta\over m}$ and $\frac{\alpha+\rho\beta}{n+\rho m}=1-{1\over q}$ together imply 
\bel{alpha constraint case1}
{\alpha\over n}~>~1-{1\over q}
\eeq	
and
\bel{beta formula <}
{\beta\over m}~<~1-{1\over q}.
\eeq	
These two strict inequalities will be used later.	
\subsection{Case 1: $k=\ell=0$}		
Recall that $a$ is supported in ${1\over 2}\Q$ and $|a(x,y)|\leq 2^{n+m}\vol\{\Q\}^{-1}$. Let $(x,y)\in \Q$. We have
\bel{case1 est1}
\begin{array}{lr}\ds
\left|\I_{\alpha\beta}^\rho a(x,y)\right|~\leq~\iint_{{1\over 2}\Q} |a(u,v)| \left({1\over |x-u|}\right)^{n-\alpha} \left[{1\over |x-u|^\rho+|y-v|}\right]^{m-\beta}dudv
\\\\ \ds~~~~~~~~~~~~~~~~~
~\lesssim~\sum_{s\leq L}\sum_{r\leq L}\iint_{2^{s-1}<|x-u|\leq 2^s, 2^{r-1}<|y-v|\leq2^r} {1\over\vol\{\Q\}}2^{s(\alpha-n)}\big[2^{r}+2^{\rho s}\big]^{\beta-m}dudv
\\\\ \ds~~~~~~~~~~~~~~~~~
~=~2^{-(n+m)L} \sum_{s\leq L}\sum_{r\leq L}\iint_{2^{s-1}<|x-u|\leq 2^s, 2^{r-1}<|y-v|\leq2^r} 2^{s(\alpha-n)}\big[2^{r}+2^{\rho s}\big]^{\beta-m}dudv.
\end{array}
\eeq
Split the above summation into two: $\S_1+\S_2$ $w.r.t$ $r\ge\rho s$ and $r<\rho s$. We find 
\bel{S_1}
\S_1~=~
2^{-(n+m)L}\sum_{s\leq \min\left\{ L/\rho,L\right\} } \sum_{\rho s\leq r\leq L}\iint_{2^{s-1}<|x-u|\leq 2^s, ~2^{r-1}<|y-v|\leq2^r} 2^{s(\alpha-n)}\big[2^{r}+2^{\rho s}\big]^{\beta-m}dudv
\eeq
and
\bel{S_2}
\S_2~=~2^{-(n+m)L}\sum_{r\leq \min\left\{L, \rho L\right\}}\sum_{{r\over \rho}<s\leq L}\iint_{2^{s-1}<|x-u|\leq 2^s,~ 2^{r-1}<|y-v|\leq2^r} 2^{s(\alpha-n)}\big[2^{r}+2^{\rho s}\big]^{\beta-m}dudv.
\eeq
Consider $\S_1$ in (\ref{S_1}). We have 
\bel{S_1 Est1}
\begin{array}{lr}\ds
\S_1 ~\lesssim~ 2^{-(n+m)L}\sum_{s\leq \min\left\{ L/\rho,L\right\} } \sum_{\rho s\leq r\leq L} 2^{s(\alpha-n)}\Big[2^r+2^{\rho s}\Big]^{\beta-m} 2^{s n}2^{rm}
\\\\ \ds~~~~
~\leq~ 2^{-(n+m)L}\sum_{s\leq \min\left\{ L/\rho,L\right\} } \sum_{\rho s\leq r\leq L} 2^{s\alpha}2^{r\beta}\qquad\hbox{\small{($0<\beta<m$)}}
\\\\ \ds~~~~
~\leq~\C_\beta~ 2^{-(n+m-\beta)L}\sum_{s\leq \min\left\{ L/\rho,L\right\} } 2^{s\alpha}.
	\end{array}
	\eeq
For $L>0$, $\min\left\{ L/\rho,L\right\}=L/ \rho$. $\S_1$ in (\ref{S_1 Est1}) is further bounded by  $\C_{\alpha~\beta} 2^{\big[(\alpha/\rho)+\beta-n-m\big]L}$. We find
\bel{Case1 Int Est1}
\begin{array}{lr}\ds
\iint_\Q \S_1^q dxdy~\leq~\C_{\alpha~\beta} ~2^{q\big[(\alpha/\rho)+\beta-n-m\big]L}2^{(n+m)L}
\\\\ \ds~~~~~~~~~~~~~~~~~~~~~~
~=~\C_{\alpha~\beta} ~2^{q\big[(\alpha/\rho)+\beta\big]L}2^{-(q-1)(n+m)L}~=~\C_{\alpha~\beta}~2^{(q-1)\big[{n\over \rho}+m\big]L}2^{-(q-1)(n+m)L}\quad \hbox{\small{$\left(\frac{\alpha+\rho\beta}{n+\rho m}=1-{1\over q}\right)$}}
\\\\ \ds~~~~~~~~~~~~~~~~~~~~~~
~=~ \C_{\alpha~\beta}~2^{(q-1)\big[{n\over \rho}-n\big]L}~\leq~\C_{\alpha~\beta~\rho~q}.
\end{array}
\eeq
On the other hand, (\ref{alpha constraint case1}) implies $q\alpha>(q-1)n$.
For $L\leq0$, $\min\left\{ L/\rho,L\right\}=L$. $\S_1$ in (\ref{S_1 Est1}) is further bounded by  $\C_{\alpha~\beta} 2^{\big[\alpha+\beta-n-m\big]L}$. We have
\bel{Case1 Int Est2}
\begin{array}{lr}\ds
\iint_\Q \S_1^q dxdy~\leq~\C_{\alpha~\beta} ~2^{q\big[\alpha+\beta-n-m\big]L}2^{(n+m)L}
\\\\ \ds~~~~~~~~~~~~~~~~~~~~
~=~\C_{\alpha~\beta} ~2^{q\alpha\big[1-{1\over \rho}\big]L}2^{q\big[(\alpha/\rho)+\beta\big]L}2^{-(q-1)(n+m)L}
\\\\ \ds~~~~~~~~~~~~~~~~~~~~
~=~\C_{\alpha~\beta}~2^{q\alpha\big[1-{1\over \rho}\big]L}2^{(q-1)\big[{n\over \rho}+m\big]L}2^{-(q-1)(n+m)L} \qquad \hbox{\small{$\left(\frac{\alpha+\rho\beta}{n+\rho m}=1-{1\over q}\right)$}}
\\\\ \ds~~~~~~~~~~~~~~~~~~~~
~=~ \C_{\alpha~\beta}~2^{q\alpha\big[1-{1\over \rho}\big]L}2^{(q-1)\big[{n\over \rho}-n\big]L}
\\\\ \ds~~~~~~~~~~~~~~~~~~~~
~=~\C_{\alpha~\beta}~2^{\big[q\alpha-(q-1)n\big]\big[1-{1\over \rho}\big]L}
\\\\ \ds~~~~~~~~~~~~~~~~~~~~
~\leq~\C_{\alpha~\beta~\rho~q}.
\end{array}
\eeq

Consider $\S_2$ in (\ref{S_2}). We have 
\bel{S_2 Est1}
\begin{array}{lr}\ds
\S_2 ~\lesssim~ 2^{-(n+m)L}\sum_{r\leq \min\left\{L, \rho L\right\}}\sum_{{r\over \rho}<s\leq L} 2^{s(\alpha-n)}\Big[2^r+2^{\rho s}\Big]^{\beta-m} 2^{s n}2^{rm}
\\\\ \ds~~~~
~\leq~ 2^{-(n+m)L}\sum_{r\leq \min\left\{L, \rho L\right\}}\sum_{{r\over \rho}<s\leq L} 2^{s\big[\alpha+\rho\beta-\rho m\big]}2^{rm}\qquad\hbox{\small{($0<\beta<m$)}}
\\\\ \ds~~~~
~\leq~\C_{\alpha~\beta~\rho}~ 2^{-(n+m)L}\sum_{r\leq \min\left\{L, \rho L\right\}} \max\Bigg\{2^{(r/\rho)(\alpha+\rho\beta)},2^{L\big[\alpha+\rho\beta-\rho m\big]}2^{rm}\Bigg\}.
	\end{array}
	\eeq
For $L>0$, $\min\left\{L, \rho L\right\}=L$. $\S_2$ in (\ref{S_2 Est1}) is further bounded by 
$\C_{\alpha~\beta~\rho}2^{-(n+m)L}\max\Bigg\{2^{\big[(\alpha/\rho)+\beta\big]L},2^{\big[\alpha+\rho\beta-\rho m\big]L}2^{mL}\Bigg\}$. We find
\bel{Case1 Int Est3}
\begin{array}{lr}\ds
\iint_\Q \S_2^q dxdy~\leq~\C_{\alpha~\beta~\rho} ~\max\Bigg\{ 2^{q\big[(\alpha/\rho)+\beta\big]L}, 2^{q\big[\alpha+\rho\beta-\rho m\big]L}2^{qmL}\Bigg\} 2^{-(q-1)(n+m)L}
\\\\ \ds~~~~~~~~~~~~~~~~~~~~
~=~\C_{\alpha~\beta~\rho}~\max\Bigg\{2^{(q-1)\big[{n\over \rho}+m\big]L}, 2^{\big[(q-1)n-\rho m\big]L} 2^{qmL}\Bigg\}2^{-(q-1)(n+m)L}\qquad \hbox{\small{$\left(\frac{\alpha+\rho\beta}{n+\rho m}=1-{1\over q}\right)$}}
\\\\ \ds~~~~~~~~~~~~~~~~~~~~
~=~ \C_{\alpha~\beta~\rho}~\max\Bigg\{2^{(q-1)\big[{n\over \rho}-n\big]L}, 2^{(1-\rho)m L}\Bigg\}
~\leq~\C_{\alpha~\beta~\rho~q}.
\end{array}
\eeq
For $L\leq0$, $\min\left\{L, \rho L\right\}=\rho L$. $\S_2$ in (\ref{S_2 Est1}) is further bounded by 
$\C_{\alpha~\beta~\rho}2^{-(n+m)L}2^{(\alpha+\rho\beta)L}$. We have
\bel{Case1 Int Est3}
\begin{array}{lr}\ds
\iint_\Q \S_2^q dxdy~\leq~\C_{\alpha~\beta~\rho}~ 2^{q(\alpha+\rho\beta)L} 2^{-(q-1)(n+m)L}
\\\\ \ds~~~~~~~~~~~~~~~~~~~~~
~=~\C_{\alpha~\beta~\rho}~2^{(q-1)(n+\rho m)L} 2^{-(q-1)(n+m)L}\qquad \hbox{\small{$\left(\frac{\alpha+\rho\beta}{n+\rho m}=1-{1\over q}\right)$}}
\\\\ \ds~~~~~~~~~~~~~~~~~~~~~
~=~\C_{\alpha~\beta~\rho}~2^{(q-1)(\rho-1)mL}~\leq~\C_{\alpha~\beta~\rho~q}.
\end{array}
\eeq

\subsection{Case 2: $k>0$, $\ell>0$}		
Let $(x,y)\in \Q_{k\ell}$. We have $|x-u|\sim 2^{k+L}$, $|y-v|\sim 2^{\ell+L}$. Consider $k+L\geq0$ or $\ell\geq k$. 
We find
\bel{Case2 ineq >}
|x-u|^\rho+|y-v|~\gtrsim~|x-u|.
\eeq
By using the cancellation property of $a$: $\ds\iint_{{1\over 2}\Q} a(x,y)dxdy=0$, we have
\bel{Case2 Est1}
\I_{\alpha\beta}^\rho a(x,y)~=~ \iint_{{1\over 2}\Q} a(u,v) \Big[\Omega^{\alpha\beta}_\rho(x-u,y-v)-\Omega^{\alpha\beta}_\rho(x,y)\Big] dudv.
\eeq
Recall the estimate in (\ref{deri Est}). From (\ref{Case2 Est1}), we find
\bel{Case2 Est2}
\begin{array}{lr}\ds
\left|\I_{\alpha\beta}^\rho a(x,y)\right|~\leq~ \iint_{{1\over 2}\Q}{1\over \vol\{\Q\}} \sqrt{|u|^2+|v|^2} \left|\nabla \Omega^{\alpha\beta}_\rho(x-u,y-v)\right|dudv
\\\\ \ds~~~~~~~~~~~~~~~~~
~\leq~\C_{\alpha~\beta~\rho} \iint_{{1\over 2}\Q}{1\over \vol\{\Q\}} \sqrt{|u|^2+|v|^2} \left({1\over |x-u|}\right)^{n-\alpha} \left[{1\over |x-u|^\rho+|y-v|}\right]^{m-\beta}
\\ \ds~~~~~~~~~~~~~~~~~~~~~~~~~~~~~~~~~~~~~~~~~~~~~~~~~~~~~~~~~~~~~~~~~~~~~~~~~
\max\left\{{1\over|x-u|},{1\over |x-u|^{\rho}+|y-v|}\right\} dudv
\\\\ \ds~~~~~~~~~~~~~~~~~
~\leq~\C_{\alpha~\beta~\rho} \iint_{{1\over 2}\Q}{1\over \vol\{\Q\}} \sqrt{|u|^2+|v|^2} \left({1\over |x-u|}\right)^{n-\alpha} \left[{1\over |x-u|^\rho+|y-v|}\right]^{m-\beta} {1\over|x-u|} dudv\qquad \hbox{\small{by (\ref{Case2 ineq >})}}
\\\\ \ds~~~~~~~~~~~~~~~~~
~\leq~\C_{\alpha~\beta~\rho} \iint_{{1\over 2}\Q} 2^{-(n+m)L}2^L 2^{(k+L)(\alpha-n)} \Big[2^{\rho(k+L)}+2^{\ell+L}\Big]^{\beta-m} 2^{-(k+L)} dudv
\\\\ \ds~~~~~~~~~~~~~~~~~
~\leq~\C_{\alpha~\beta~\rho} ~2^{-k} 2^{(k+L)(\alpha-n)} \Big[2^{\rho(k+L)}+2^{\ell+L}\Big]^{\beta-m}. 
 \end{array}
\eeq
Suppose $\rho(k+L)\geq \ell+L$. The last line of (\ref{Case2 Est2}) can be further bounded by $\C_{\alpha~\beta~\rho} ~2^{-k} 2^{(k+L)(\alpha-n)} 2^{\rho(k+L)(\beta-m)}$. We have
\bel{Case2 Est3}
\begin{array}{lr}\ds
\iint_{\Q_{k\ell}} \left|\I_{\alpha\beta}^\rho a(x,y)\right|^q dxdy
\\\\ \ds
~\leq~
\C_{\alpha~\beta~\rho} \iint_{\Q_{k\ell}} 2^{-qk} 2^{q(k+L)(\alpha-n)} 2^{q\rho(k+L)(\beta-m)} dxdy
\\\\ \ds
~=~\C_{\alpha~\beta~\rho} \iint_{\Q_{k\ell}} 2^{-qk} 2^{q(k+L)\big[\alpha+\rho\beta-(n+\rho m)\big]} dxdy
\\\\ \ds
~\leq~\C_{\alpha~\beta~\rho}~2^{-qk} 2^{q(k+L)\big[\alpha+\rho\beta-(n+\rho m)\big]} 2^{(k+L)n}2^{(\ell+L)m}.
\end{array}
\eeq
By using (\ref{Case2 Est3}) and taking the summation over every $k,\ell\ge0$: 
$\rho(k+L)\geq \ell+L$, we obtain 
\bel{Case2 Est4}
\begin{array}{lr}\ds
\sum_{k,\ell\ge0\colon \rho(k+L)\geq \ell+L} \iint_{\Q_{k\ell}} \left|\I_{\alpha\beta}^\rho a(x,y)\right|^q dxdy
\\\\ \ds
~\leq~\C_{\alpha~\beta~\rho}\sum_{k\ge0}\sum_{\ell\leq\rho(k+L)-L} 2^{-qk} 2^{q(k+L)\big[\alpha+\rho\beta-(n+\rho m)\big]} 2^{(k+L)n}2^{(\ell+L)m}
\\\\ \ds
~\leq~ \C_{\alpha~\beta~\rho}\sum_{k\ge0}2^{-qk} 2^{q(k+L)\big[\alpha+\rho\beta-(n+\rho m)\big]} 2^{(k+L)n}2^{\rho(k+L)m}
\\\\ \ds
~=~\C_{\alpha~\beta~\rho}\sum_{k\ge0}2^{-qk}  2^{(k+L) \big[q(\alpha+\rho\beta)-(q-1)(n+\rho m)\big]}
\\\\ \ds
~=~\C_{\alpha~\beta~\rho}\sum_{k\ge0}2^{-qk} \qquad\hbox{\small{$\left({\alpha+\rho\beta\over n+\rho m}=1-{1\over q}\right)$}}
\\\\ \ds
~\leq~\C_{\alpha~\beta~\rho~q}.
\end{array}
\eeq
Suppose $\rho(k+L)\leq \ell+L$. The last line of (\ref{Case2 Est2}) is further bounded by $\C_{\alpha~\beta~\rho} ~2^{-k} 2^{(k+L)(\alpha-n)} 2^{(\ell+L)(\beta-m)}$. We have		
\bel{Case2 Est5}
\begin{array}{lr}\ds
\iint_{\Q_{k\ell}} \left|\I_{\alpha\beta}^\rho a(x,y)\right|^q dxdy
\\\\ \ds
~\leq~
\C_{\alpha~\beta~\rho} \iint_{\Q_{k\ell}} 2^{-qk} 2^{q(k+L)(\alpha-n)} 2^{q(\ell+L)(\beta-m)} dxdy
\\\\ \ds
~\leq~\C_{\alpha~\beta~\rho}~2^{-qk} 2^{q(k+L)(\alpha-n)} 2^{q(\ell+L)(\beta-m)} 2^{(k+L)n}2^{(\ell+L)m}
\\\\ \ds
~=~\C_{\alpha~\beta~\rho}~2^{-qk} 2^{q(k+L)(\alpha-n)} 2^{(\ell+L)\big[q(\beta-m)+m\big]} 2^{(k+L)n}.
\end{array}
\eeq
Recall (\ref{beta formula <}). 
Note that ${\beta\over m}<1-{1\over q}$ implies $q(\beta-m)<-m$.

By using (\ref{Case2 Est5}) and taking the summation over every $k,\ell\ge0$: 
$\rho(k+L)\leq \ell+L$, we obtain 
\bel{Case2 Est6}
\begin{array}{lr}\ds
\sum_{k,\ell\ge0\colon \rho(k+L)\leq \ell+L} \iint_{\Q_{k\ell}} \left|\I_{\alpha\beta}^\rho a(x,y)\right|^q dxdy
\\\\ \ds
~\leq~\C_{\alpha~\beta~\rho}\sum_{k\ge0}\sum_{\ell\ge\rho(k+L)-L} 2^{-qk} 2^{q(k+L)(\alpha-n)} 2^{(\ell+L)\big[q(\beta-m)+m\big]} 2^{(k+L)n}
\\\\ \ds
~\leq~ \C_{\alpha~\beta~\rho}\sum_{k\ge0}2^{-qk} 2^{q(k+L)(\alpha-n)} 2^{\rho(k+L)\big[q(\beta-m)+m\big]} 2^{(k+L)n}\qquad\hbox{\small{($q(\beta-m)<-m$)}}
\\\\ \ds
~=~\C_{\alpha~\beta~\rho}\sum_{k\ge0}2^{-qk}  2^{(k+L) \big[q(\alpha+\rho\beta)-(q-1)(n+\rho m)\big]}
\\\\ \ds
~=~\C_{\alpha~\beta~\rho}\sum_{k\ge0}2^{-qk}~\leq~\C_{\alpha~\beta~\rho~q}.\qquad\hbox{\small{$\left({\alpha+\rho\beta\over n+\rho m}=1-{1\over q}\right)$}}
\end{array}
\eeq
Note that $|x-u|\sim 2^{k+L}$, $|y-v|\sim 2^{\ell+L}$ for $(x,y)\in \Q_{k\ell}$.
Consider $k+L<0$ and $\ell<k$. We find
\bel{Case2 ineq <}
|x-u|^\rho+|y-v|~\lesssim~  |x-u|.
\eeq
Recall  (\ref{deri Est}) and (\ref{Case2 Est1}). We have
\bel{Case2 Est7}
\begin{array}{lr}\ds
\left|\I_{\alpha\beta}^\rho a(x,y)\right|~\leq~ \iint_{{1\over 2}\Q}{1\over \vol\{\Q\}} \sqrt{|u|^2+|v|^2} \left|\nabla \Omega^{\alpha\beta}_\rho(x-u,y-v)\right|dudv
\\\\ \ds~~~~~~~~~~~~~~~~~
~\leq~\C_{\alpha~\beta~\rho} \iint_{{1\over 2}\Q}{1\over \vol\{\Q\}} \sqrt{|u|^2+|v|^2} \left({1\over |x-u|}\right)^{n-\alpha} \left[{1\over |x-u|^\rho+|y-v|}\right]^{m-\beta}
\\ \ds~~~~~~~~~~~~~~~~~~~~~~~~~~~~~~~~~~~~~~~~~~~~~~~~~~~~~~~~~~~~~~~~~~~~~~~~~
\max\left\{{1\over|x-u|},{1\over |x-u|^{\rho}+|y-v|}\right\} dudv
\\\\ \ds~~~~~~~~~~~~~~~~~
~\leq~\C_{\alpha~\beta~\rho} \iint_{{1\over 2}\Q}{1\over \vol\{\Q\}} \sqrt{|u|^2+|v|^2} \left({1\over |x-u|}\right)^{n-\alpha} \left[{1\over |x-u|^\rho+|y-v|}\right]^{m-\beta} {1\over|x-u|^\rho+|y-v|} dudv
\\ \ds~~~~~~~~~~~~~~~~~~~~~~~~~~~~~~~~~~~~~~~~~~~~~~~~~~~~~~~~~~~~~~~~~~~~~~~~~~~~~~~~~~~~~~~~~~~~~~~~~~~~~~~~~~~~~~~~~~~~~~~~~~~~~~~~~~~~~~~~~~~~~~~~~~~~~~~~~~~~~
 \hbox{\small{by (\ref{Case2 ineq <})}}
\\ \ds~~~~~~~~~~~~~~~~~
~\leq~\C_{\alpha~\beta~\rho} \iint_{{1\over 2}\Q} 2^{-(n+m)L}2^L 2^{(k+L)(\alpha-n)} \Big[2^{\rho(k+L)}+2^{\ell+L}\Big]^{\beta-m} \Big[2^{\rho(k+L)}+2^{\ell+L}\Big]^{-1} dudv
\\\\ \ds~~~~~~~~~~~~~~~~~
~\leq~\C_{\alpha~\beta~\rho} ~2^L 2^{(k+L)(\alpha-n)} \Big[2^{\rho(k+L)}+2^{\ell+L}\Big]^{\beta-m} \Big[2^{\rho(k+L)}+2^{\ell+L}\Big]^{-1}. 
 \end{array}
\eeq
Suppose $\rho(k+L)\geq \ell+L$. The last line in (\ref{Case2 Est7}) is further bounded by $\C_{\alpha~\beta~\rho} 2^L 2^{(k+L)(\alpha-n)} 2^{\rho(k+L)(\beta-m)} 2^{-\rho(k+L)}$.
We have
\bel{Case2 Est8}
\begin{array}{lr}\ds
\iint_{\Q_{k\ell}}\left|\I_{\alpha\beta}^\rho a (x,y)\right|^qdxdy ~\leq~ \C_{\alpha~\beta~\rho} \iint_{\Q_{k\ell}}
2^{qL} 2^{q(k+L)(\alpha-n)} 2^{q\rho(k+L)(\beta-m)} 2^{-q\rho(k+L)} dxdy
\\\\ \ds~~~~~~~~~~~~~~~~~~~~~~~~~~~~~~~~~~~~~
~\leq~\C_{\alpha~\beta~\rho} ~
2^{qL} 2^{q(k+L)(\alpha-n)} 2^{q\rho(k+L)(\beta-m)} 2^{-q\rho(k+L)} 2^{(k+L)n}2^{(\ell+L)m}
\\\\ \ds~~~~~~~~~~~~~~~~~~~~~~~~~~~~~~~~~~~~~
~=~\C_{\alpha~\beta~\rho} ~2^{q\big[L-\rho(k+L)\big]} 2^{(k+L)\big[ q(\alpha+\rho\beta-n-\rho m)+n\big]} 2^{(\ell+L)m}.
\end{array}
\eeq
By using (\ref{Case2 Est8}) and summing over every $k,\ell\ge0$: $\rho(k+L)\geq \ell+L$, we obtain
\bel{Case2 Est9}
\begin{array}{lr}\ds
 \sum_{k,\ell\ge0\colon \rho(k+L)\geq \ell+L} \iint_{\Q_{k\ell}}\left|\I_{\alpha\beta}^\rho a (x,y)\right|^qdxdy
 \\\\ \ds
 ~\leq~\C_{\alpha~\beta~\rho} \sum_{k\ge-(\rho-1)L/\rho} \sum_{\ell\leq \rho(k+L)-L} 2^{q\big[L-\rho(k+L)\big]} 2^{(k+L)\big[ q(\alpha+\rho\beta-n-\rho m)+n\big]} 2^{(\ell+L)m} 
  \\\\ \ds
 ~\leq~\C_{\alpha~\beta~\rho} \sum_{k\ge-(\rho-1)L/\rho}  2^{q\big[L-\rho(k+L)\big]} 2^{(k+L)\big[ q(\alpha+\rho\beta-n-\rho m)+n\big]} 2^{\rho(k+L)m} 
\\\\ \ds
~=~\C_{\alpha~\beta~\rho} \sum_{k\ge-(\rho-1)L/\rho}  2^{q\big[L-\rho(k+L)\big]}2^{(k+L)\big[q(\alpha+\rho m)-(q-1)(n+\rho m)\big]}
\\\\ \ds
~=~ \C_{\alpha~\beta~\rho} \sum_{k\ge-(\rho-1)L/\rho} 2^{q\big[L-\rho(k+L)\big]}~\leq~\C_{\alpha~\beta~\rho~q}.\qquad \hbox{\small{$\left({\alpha+\rho\beta\over n+\rho m}=1-{1\over q}\right)$}}			
\end{array}
\eeq
Suppose $\rho(k+L)<\ell+L$. The last line in (\ref{Case2 Est7}) is further bounded by $\C_{\alpha~\beta~\rho} ~2^L 2^{(k+L)(\alpha-n)} 2^{(\ell+L)(\beta-m)}2^{-(\ell+L)}$. 
\bel{Case2 Est10}
\begin{array}{lr}\ds
\iint_{\Q_{k\ell}}\left|\I_{\alpha\beta}^\rho a (x,y)\right|^qdxdy ~\leq~ \C_{\alpha~\beta~\rho} \iint_{\Q_{k\ell}}
2^{qL} 2^{q(k+L)(\alpha-n)} 2^{q(\ell+L)(\beta-m)}2^{-q(\ell+L)} dxdy
\\\\ \ds~~~~~~~~~~~~~~~~~~~~~~~~~~~~~~~~~~~~~~~
~\leq~\C_{\alpha~\beta~\rho} 
2^{qL} 2^{q(k+L)(\alpha-n)} 2^{q(\ell+L)(\beta-m)}2^{-q(\ell+L)} 2^{(k+L)n}2^{(\ell+L)m}
\\\\ \ds~~~~~~~~~~~~~~~~~~~~~~~~~~~~~~~~~~~~~~~
~=~\C_{\alpha~\beta~\rho} ~2^{-q\ell} 2^{(k+L)\big[q(\alpha-n)+n\big]} 2^{(\ell+L)\big[q(\beta-m)+m\big]}.
\end{array}
\eeq
Recall (\ref{alpha constraint case1}). Note that ${\alpha\over n}>1-{1\over q}$ implies $q(\alpha-n)>-n$. 
By using (\ref{Case2 Est10}) and summing over every $k,\ell\ge0$: $\rho(k+L)< \ell+L$, we obtain
\bel{Case2 Est11}
\begin{array}{lr}\ds
 \sum_{k,\ell\ge0\colon \rho(k+L)< \ell+L} \iint_{\Q_{k\ell}}\left|\I_{\alpha\beta}^\rho a (x,y)\right|^qdxdy
 \\\\ \ds
 ~\leq~\C_{\alpha~\beta~\rho}  \sum_{\ell\ge0} \sum_{k<\ell/\rho+L/\rho-L} 2^{-q\ell} 2^{(k+L)\big[q(\alpha-n)+n\big]} 2^{(\ell+L)\big[q(\beta-m)+m\big]}
 
   \\\\ \ds
 ~\leq~\C_{\alpha~\beta~\rho} \sum_{\ell\ge0} 2^{-q\ell} 2^{\big[\ell/\rho+L/\rho\big]\big[q(\alpha-n)+n\big]} 2^{(\ell+L)\big[q(\beta-m)+m\big]}\qquad\hbox{\small{($q(\alpha-n)>-n$)}}
\\\\ \ds
~=~\C_{\alpha~\beta~\rho} \sum_{\ell\ge0}  2^{-q\ell} 2^{(\ell+L)\big[q(\alpha+\rho\beta)-(q-1)(n+\rho m)\big]/\rho}
\\\\ \ds
~=~ \C_{\alpha~\beta~\rho} \sum_{\ell\ge0} 2^{-q\ell}~\leq~\C_{\alpha~\beta~\rho~q}.\qquad \hbox{\small{$\left({\alpha+\rho\beta\over n+\rho m}=1-{1\over q}\right)$}}			
\end{array}
\eeq

\subsection{Case 3: $k>0$, $\ell=0$ or $k=0$, $\ell>0$}		
Let $(x,y)\in \Q_{k0}$. We have $|x-u|\sim 2^{k+L}$. Consider $\rho(k+L)\geq L$. If $k+L\geq0$, 
we find
\bel{Case3 ineq >}
|x-u|^\rho~\gtrsim~|x-u|.
\eeq
Recall  (\ref{deri Est}) and (\ref{Case2 Est1}). We have
\bel{Case3 Est1}
\begin{array}{lr}\ds
\left|\I_{\alpha\beta}^\rho a(x,y)\right|~\leq~ \iint_{{1\over 2}\Q}{1\over \vol\{\Q\}} \sqrt{|u|^2+|v|^2} \left|\nabla \Omega^{\alpha\beta}_\rho(x-u,y-v)\right|dudv
\\\\ \ds
\leq~\C_{\alpha~\beta~\rho} \iint_{\Q}{1\over \vol\{{1\over 2}\Q\}} \sqrt{|u|^2+|v|^2} \left({1\over |x-u|}\right)^{n-\alpha} \left[{1\over |x-u|^\rho+|y-v|}\right]^{m-\beta}
\max\left\{{1\over|x-u|},{1\over |x-u|^{\rho}+|y-v|}\right\} dudv
\\\\ \ds
\leq~\C_{\alpha~\beta~\rho} \iint_{{1\over 2}\Q}{1\over \vol\{\Q\}} \sqrt{|u|^2+|v|^2} \left({1\over |x-u|}\right)^{n-\alpha} \left[{1\over |x-u|^\rho+|y-v|}\right]^{m-\beta} {1\over|x-u|} dudv\qquad \hbox{\small{by (\ref{Case3 ineq >})}}
\\\\ \ds
\leq~\C_{\alpha~\beta~\rho} \iint_{{1\over 2}\Q} 2^{-(n+m)L}2^L 2^{(k+L)(\alpha-n)} \Big[2^{\rho(k+L)}+2^L\Big]^{\beta-m} 2^{-(k+L)} dudv
\\\\ \ds
\leq~\C_{\alpha~\beta~\rho} ~2^{-k} 2^{(k+L)(\alpha-n)} \Big[2^{\rho(k+L)}+2^L\Big]^{\beta-m}. 
 \end{array}
\eeq
 The last line of (\ref{Case3 Est1}) is further bounded by $\C_{\alpha~\beta~\rho} ~2^{-k} 2^{(k+L)(\alpha-n)} 2^{\rho(k+L)(\beta-m)}$. We have
\bel{Case3 Est2}
\begin{array}{lr}\ds
\iint_{\Q_{k0}} \left|\I_{\alpha\beta}^\rho a(x,y)\right|^q dxdy~\leq~
\C_{\alpha~\beta~\rho} \iint_{\Q_{k0}} 2^{-qk} 2^{q(k+L)(\alpha-n)} 2^{q\rho(k+L)(\beta-m)} dxdy
\\\\ \ds~~~~~~~~~~~~~~~~~~~~~~~~~~~~~~~~~~~~~
~=~\C_{\alpha~\beta~\rho} \iint_{\Q_{k0}} 2^{-qk} 2^{q(k+L)\big[\alpha+\rho\beta-(n+\rho m)\big]} dxdy
\\\\ \ds~~~~~~~~~~~~~~~~~~~~~~~~~~~~~~~~~~~~~
~\leq~\C_{\alpha~\beta~\rho}~2^{-qk} 2^{q(k+L)\big[\alpha+\rho\beta-(n+\rho m)\big]} 2^{(k+L)n}2^{Lm}.
\end{array}
\eeq
By using (\ref{Case3 Est2}) and taking the summation over every $k\ge0$: 
$\rho(k+L)\geq L$, we obtain 
\bel{Case2 Est3}
\begin{array}{lr}\ds
\sum_{k\ge0\colon \rho(k+L)\geq L} \iint_{\Q_{k0}} \left|\I_{\alpha\beta}^\rho a(x,y)\right|^q dxdy
\\\\ \ds
~\leq~\C_{\alpha~\beta~\rho} \sum_{k\ge0\colon \rho(k+L)\geq L} 2^{-qk} 2^{q(k+L)\big[\alpha+\rho\beta-(n+\rho m)\big]} 2^{(k+L)n} 2^{Lm}
\\\\ \ds
~\leq~ \C_{\alpha~\beta~\rho}\sum_{k\ge0}2^{-qk} 2^{q(k+L)\big[\alpha+\rho\beta-(n+\rho m)\big]} 2^{(k+L)n}2^{\rho(k+L)m}
\\\\ \ds
~=~\C_{\alpha~\beta~\rho}\sum_{k\ge0}2^{-qk}  2^{(k+L) \big[q(\alpha+\rho\beta)-(q-1)(n+\rho m)\big]}
\\\\ \ds
~=~\C_{\alpha~\beta~\rho}\sum_{k\ge0}2^{-qk}~\leq~\C_{\alpha~\beta~\rho~q}.\qquad\hbox{\small{$\left({\alpha+\rho\beta\over n+\rho m}=1-{1\over q}\right)$}}
\end{array}
\eeq
On the other hand, if $k+L<0$, we find
\bel{Case3 ineq <}
|x-u|^\rho+|y-v|~\lesssim~2^{\rho(k+L)}+2^L~\lesssim~2^{k+L}~\approx~|x-u|.
\eeq
Recall  (\ref{deri Est}) and (\ref{Case2 Est1}). We have
\bel{Case3 Est4}
\begin{array}{lr}\ds
\left|\I_{\alpha\beta}^\rho a(x,y)\right|~\leq~ \iint_{{1\over 2}\Q}{1\over \vol\{\Q\}} \sqrt{|u|^2+|v|^2} \left|\nabla \Omega^{\alpha\beta}_\rho(x-u,y-v)\right|dudv
\\\\ \ds~~~~~~~~~~~~~~~~~
~\leq~\C_{\alpha~\beta~\rho} \iint_{{1\over 2}\Q}{1\over \vol\{\Q\}} \sqrt{|u|^2+|v|^2} \left({1\over |x-u|}\right)^{n-\alpha} \left[{1\over |x-u|^\rho+|y-v|}\right]^{m-\beta}
\\ \ds~~~~~~~~~~~~~~~~~~~~~~~~~~~~~~~~~~~~~~~~~~~~~~~~~~~~~~~~~~~~~~~~~~~~~~~~~
\max\left\{{1\over|x-u|},{1\over |x-u|^{\rho}+|y-v|}\right\} dudv
\\\\ \ds~~~~~~~~~~~~~~~~~
~\leq~\C_{\alpha~\beta~\rho} \iint_{{1\over 2}\Q}{1\over \vol\{\Q\}} \sqrt{|u|^2+|v|^2} \left({1\over |x-u|}\right)^{n-\alpha} \left[{1\over |x-u|^\rho+|y-v|}\right]^{m-\beta} {1\over|x-u|^\rho+|y-v|} dudv
\\ \ds~~~~~~~~~~~~~~~~~~~~~~~~~~~~~~~~~~~~~~~~~~~~~~~~~~~~~~~~~~~~~~~~~~~~~~~~~~~~~~~~~~~~~~~~~~~~~~~~~~~~~~~~~~~~~~~~~~~~~~~~~~~~~~~~~~~~~~~~~~~~~~~~~~~~~~~~~~~~
 \hbox{\small{by (\ref{Case3 ineq <})}}
\\ \ds~~~~~~~~~~~~~~~~~
~\leq~\C_{\alpha~\beta~\rho} \iint_{{1\over 2}\Q} 2^{-(n+m)L}2^L 2^{(k+L)(\alpha-n)} \Big[2^{\rho(k+L)}+2^{L}\Big]^{\beta-m} \Big[2^{\rho(k+L)}+2^{L}\Big]^{-1} dudv
\\\\ \ds~~~~~~~~~~~~~~~~~
~\leq~\C_{\alpha~\beta~\rho} ~2^L 2^{(k+L)(\alpha-n)} \Big[2^{\rho(k+L)}+2^{L}\Big]^{\beta-m} \Big[2^{\rho(k+L)}+2^{L}\Big]^{-1}. 
 \end{array}
\eeq
Because $\rho(k+L)\geq L$, the last line in (\ref{Case3 Est4}) can be further bounded by $\C_{\alpha~\beta~\rho} 2^L 2^{(k+L)(\alpha-n)} 2^{\rho(k+L)(\beta-m)} 2^{-\rho(k+L)}$.
We have
\bel{Case3 Est5}
\begin{array}{lr}\ds
\iint_{\Q_{k\ell}}\left|\I_{\alpha\beta}^\rho a (x,y)\right|^qdxdy ~\leq~ \C_{\alpha~\beta~\rho} \iint_{\Q_{k\ell}}
2^{qL} 2^{q(k+L)(\alpha-n)} 2^{q\rho(k+L)(\beta-m)} 2^{-q\rho(k+L)} dxdy
\\\\ \ds~~~~~~~~~~~~~~~~~~~~~~~~~~~~~~~~~~~~~~
~\leq~\C_{\alpha~\beta~\rho} 
2^{qL} 2^{q(k+L)(\alpha-n)} 2^{q\rho(k+L)(\beta-m)} 2^{-q\rho(k+L)} 2^{(k+L)n}2^{Lm}
\\\\ \ds~~~~~~~~~~~~~~~~~~~~~~~~~~~~~~~~~~~~~~
~=~\C_{\alpha~\beta~\rho} ~2^{q\big[L-\rho(k+L)\big]} 2^{(k+L)\big[ q(\alpha+\rho\beta-n-\rho m)+n\big]} 2^{Lm}.
\end{array}
\eeq
By using (\ref{Case3 Est5}) and summing over every $k\ge0$: $\rho(k+L)\geq L$, we obtain
\bel{Case3 Est6}
\begin{array}{lr}\ds
 \sum_{k\ge0\colon \rho(k+L)\geq L} \iint_{\Q_{k0}}\left|\I_{\alpha\beta}^\rho a (x,y)\right|^qdxdy
 \\\\ \ds
 ~\leq~\C_{\alpha~\beta~\rho} \sum_{k\ge0\colon \rho(k+L)\geq L}  2^{q\big[L-\rho(k+L)\big]} 2^{(k+L)\big[ q(\alpha+\rho\beta-n-\rho m)+n\big]} 2^{Lm} 
  \\\\ \ds
 ~\leq~\C_{\alpha~\beta~\rho} \sum_{k\ge-(\rho-1)L/\rho}  2^{q\big[L-\rho(k+L)\big]} 2^{(k+L)\big[ q(\alpha+\rho\beta-n-\rho m)+n\big]} 2^{\rho(k+L)m} 
\\\\ \ds
~=~\C_{\alpha~\beta~\rho} \sum_{k\ge-(\rho-1)L/\rho}  2^{q\big[L-\rho(k+L)\big]}2^{(k+L)\big[q(\alpha+\rho m)-(q-1)(n+\rho m)\big]}
\\\\ \ds
~=~ \C_{\alpha~\beta~\rho} \sum_{k\ge-(\rho-1)L/\rho} 2^{q\big[L-\rho(k+L)\big]}~\leq~\C_{\alpha~\beta~\rho~q}.\qquad \hbox{\small{$\left({\alpha+\rho\beta\over n+\rho m}=1-{1\over q}\right)$}}			
\end{array}
\eeq

Consider $\rho(k+L)<L$. Recall $\I_{\alpha\beta}^\rho f=f\ast\Omega^{\alpha\beta}_\rho$ defined in (\ref{I abrho})-(\ref{Omega}). We have
\bel{Case3 Est7}
\begin{array}{lr}\ds	
\left|\I_{\alpha\beta}^\rho a(x,y)\right|~\leq~ \iint_{{1\over 2}\Q}{1\over \vol\{\Q\}}   \Omega^{\alpha\beta}_\rho(x-u,y-v)dudv
\\\\ \ds~~~~~~~~~~~~~~~~~
~\leq~\sum_{s\leq L} \iint_{{1\over 2}\Q\cap\{\sqrt{m}2^{s-1}\leq|y-v|<\sqrt{m}2^s\}}{1\over \vol\{\Q\}}   \Omega^{\alpha\beta}_\rho(x-u,y-v)dudv
\\\\ \ds~~~~~~~~~~~~~~~~~
~\leq~\C_{\alpha~\beta~\rho} \sum_{s\leq L} 2^{-(n+m)L} 2^{(k+L)(\alpha-n)}\Big[2^{\rho(k+L)}+2^s\Big]^{\beta-m} 2^{Ln}2^{sm}
\\\\ \ds~~~~~~~~~~~~~~~~~
~=~\C_{\alpha~\beta~\rho} \sum_{s\leq L} 2^{-Lm} 2^{(k+L)(\alpha-n)}\Big[2^{\rho(k+L)}+2^s\Big]^{\beta-m} 2^{sm}
\\\\ \ds~~~~~~~~~~~~~~~~~
~\leq~\C_{\alpha~\beta~\rho} \sum_{s\leq L} 2^{-Lm} 2^{(k+L)(\alpha-n)} 2^{s\beta}
\\\\ \ds~~~~~~~~~~~~~~~~~
~\leq~\C_{\alpha~\beta~\rho} ~2^{(k+L)(\alpha-n)}2^{L(\beta-m)}.
\end{array}
\eeq
By using (\ref{Case3 Est7}), we find
\bel{Case3 Est8}
\begin{array}{lr}\ds	
\iint_{\Q_{k0}} 	\left|\I_{\alpha\beta}^\rho a(x,y)\right|^q dxdy~\leq~	\C_{\alpha~\beta~\rho} \iint_{\Q_{k0}}2^{q(k+L)(\alpha-n)}2^{qL(\beta-m)}dxdy
\\\\ \ds~~~~~~~~~~~~~~~~~~~~~~~~~~~~~~~~~~~~~~
~\leq~\C_{\alpha~\beta~\rho}~ 2^{q(k+L)(\alpha-n)}2^{qL(\beta-m)} 2^{(k+L)n}2^{Lm}
\\\\ \ds~~~~~~~~~~~~~~~~~~~~~~~~~~~~~~~~~~~~~~
~=~\C_{\alpha~\beta~\rho}~2^{(k+L)\big[q(\alpha-n)+n\big]} 2^{L\big[q(\beta-m)+m\big]}.
\end{array}
\eeq
Recall (\ref{alpha constraint case1}). Note that ${\alpha\over n}>1-{1\over q}$ implies $q(\alpha-n)>-n$.	
By using (\ref{Case3 Est8}) and summation over $k\ge0\colon \rho(k+L)<L$, we obtain 
\bel{Case3 Estq}
\begin{array}{lr}\ds	
\sum_{k\ge0\colon \rho(k+L)<L}\iint_{\Q_{k0}} 	\left|\I_{\alpha\beta}^\rho a(x,y)\right|^q dxdy
~\leq~	\C_{\alpha~\beta~\rho} \sum_{k\ge0\colon \rho(k+L)<L} 2^{(k+L)\big[q(\alpha-n)+n\big]} 2^{L\big[q(\beta-m)+m\big]}
\\\\ \ds~~~~~~~~~~~~~~~~~~~~~~~~~~~~~~~~~~~~~~~~~~~~~~~~~~~~~~~~~
~=~\C_{\alpha~\beta~\rho} \sum_{0\leq k<L/\rho-L} 2^{(k+L)\big[q(\alpha-n)+n\big]} 2^{L\big[q(\beta-m)+m\big]}
\\\\ \ds~~~~~~~~~~~~~~~~~~~~~~~~~~~~~~~~~~~~~~~~~~~~~~~~~~~~~~~~~
~\leq~\C_{\alpha~\beta~\rho~q}~ 2^{L\big[q(\alpha-n)+n\big]/\rho} 2^{L\big[q(\beta-m)+m\big]}\qquad\hbox{\small{($q(\alpha-n)+n>0$)}}
\\\\ \ds~~~~~~~~~~~~~~~~~~~~~~~~~~~~~~~~~~~~~~~~~~~~~~~~~~~~~~~~~
~=~\C_{\alpha~\beta~\rho~q}~2^{L\big[q(\alpha+\rho\beta)-(q-1)(n+\rho m)\big]/\rho}
\\\\ \ds~~~~~~~~~~~~~~~~~~~~~~~~~~~~~~~~~~~~~~~~~~~~~~~~~~~~~~~~~
~=~\C_{\alpha~\beta~\rho~q}.\qquad\hbox{\small{$\left( {\alpha+\rho\beta\over n+\rho m}=1-{1\over q}\right)$}}
\end{array}
\eeq

\subsection{Case 4: $k=0$ and $\ell>0$}
Let $(x,y)\in \Q_{0\ell}$. We have $|y-v|\sim 2^{\ell+L}$. Recall $\Omega^{\alpha\beta}_\rho(x,y)$ defined in (\ref{Omega}). Moreover, $\supp a\subset{1\over 2}\Q$ and $|a(x,y)|\leq2^{n+m}\vol\{\Q\}^{-1}$.	
We have
\bel{Case4 Est1}
\begin{array}{lr}\ds	
\left|\I_{\alpha\beta}^\rho a(x,y)\right| ~\leq~ \iint_{{1\over 2}\Q}|a(u,v)|   \Omega^{\alpha\beta}_\rho(x-u,y-v)dudv
\\\\ \ds~~~~~~~~~~~~~~~~~
~\lesssim~ \iint_{{1\over 2}\Q}{1\over \vol\{\Q\}}   \left({1\over |x-u|}\right)^{n-\alpha} \left[{1\over |x-u|^\rho+|y-v|}\right]^{m-\beta}dudv
\\\\ \ds~~~~~~~~~~~~~~~~~
~\leq~\sum_{s\leq L} \iint_{{1\over 2}\Q\cap\{\sqrt{n}2^{s-1}\leq|x-u|<\sqrt{n}2^s\}}{1\over \vol\{\Q\}}   \left({1\over |x-u|}\right)^{n-\alpha} \left[{1\over |x-u|^\rho+|y-v|}\right]^{m-\beta}dudv
\\\\ \ds~~~~~~~~~~~~~~~~~
~\leq~\C_{\alpha~\beta~\rho} \sum_{s\leq L} 2^{-(n+m)L} 2^{s(\alpha-n)}\Big[2^{\rho s}+2^{\ell+L}\Big]^{\beta-m} 
2^{sn}2^{Lm}
\\\\ \ds~~~~~~~~~~~~~~~~~
~=~\C_{\alpha~\beta~\rho} \sum_{s\leq L} 2^{-Ln} 2^{s\alpha}\Big[2^{\rho s}+2^{\ell+L}\Big]^{\beta-m}. 
\end{array}
\eeq
Suppose $\ell\geq (\rho-1)L$.  We find $\rho s\leq \rho L\leq \ell+L$. The last line of (\ref{Case4 Est1}) can be further bounded by
\bel{Case4 Est2}
\begin{array}{lr}\ds
 \C_{\alpha~\beta~\rho}\sum_{s\leq L} 2^{-Ln}2^{s\alpha}2^{(\ell+L)(\beta-m)}
~\leq~ \C_{\alpha~\beta~\rho}~ 2^{L(\alpha-n)}2^{(\ell+L)(\beta-m)}.
\end{array}
\eeq
From (\ref{Case4 Est1})-(\ref{Case4 Est2}), we have
\bel{Case4 Est3}
\begin{array}{lr}\ds
\iint_{\Q_{0\ell}} \left|\I_{\alpha\beta}^\rho a(x,y)\right|^q dxdy~\leq~ \C_{\alpha~\beta~\rho}\iint_{\Q_{0\ell}} 2^{qL(\alpha-n)}2^{q(\ell+L)(\beta-m)} dxdy
\\\\ \ds~~~~~~~~~~~~~~~~~~~~~~~~~~~~~~~~~~~~~~
~\leq~\C_{\alpha~\beta~\rho} ~2^{qL(\alpha-n)}2^{q(\ell+L)(\beta-m)} 2^{Ln} 2^{(\ell+L)m}
\\\\ \ds~~~~~~~~~~~~~~~~~~~~~~~~~~~~~~~~~~~~~~
~=~\C_{\alpha~\beta~\rho}~2^{L\big[ q(\alpha-n)+n\big]} 2^{(\ell+L)\big[q(\beta-m)+m\big]}.
\end{array}
\eeq
Recall (\ref{beta formula <}). We have ${\beta\over m}<1-{1\over q}\Longrightarrow q(\beta-m)+m<0$. 
By using (\ref{Case4 Est3}) and summing over $\ell\ge0\colon \ell\geq (\rho-1)L$, we obtain
\bel{Case4 Est4}
\begin{array}{lr}\ds
\sum_{\ell\ge0\colon \ell\geq (\rho-1)L} \iint_{\Q_{0\ell}} \left|\I_{\alpha\beta}^\rho a(x,y)\right|^q dxdy
~\leq~\C_{\alpha~\beta~\rho} \sum_{\ell\ge0\colon \ell\geq (\rho-1)L} 2^{L\big[ q(\alpha-n)+n\big]} 2^{(\ell+L)\big[q(\beta-m)+m\big]}
\\\\ \ds~~~~~~~~~~~~~~~~~~~~~~~~~~~~~~~~~~~~~~~~~~~~~~~~~~~~~~~~~
~\leq~\C_{\alpha~\beta~\rho~q} ~2^{L\big[ q(\alpha-n)+n\big]} 2^{\rho L\big[q(\beta-m)+m\big]}
\\\\ \ds~~~~~~~~~~~~~~~~~~~~~~~~~~~~~~~~~~~~~~~~~~~~~~~~~~~~~~~~~
~=~\C_{\alpha~\beta~\rho~q}~2^{L\big[q(\alpha+\rho\beta)-(q-1)(n+\rho m)\big]}
\\\\ \ds~~~~~~~~~~~~~~~~~~~~~~~~~~~~~~~~~~~~~~~~~~~~~~~~~~~~~~~~~
~=~\C_{\alpha~\beta~\rho~q}.\qquad \hbox{\small{$\left({\alpha+\rho\beta\over n+\rho m}=1-{1\over q}\right)$}}
\end{array}
\eeq

Suppose $\ell<(\rho-1)L$. The last line of (\ref{Case4 Est1}) can be further bounded by
\bel{Case4 Est5}
\begin{array}{lr}\ds
\C_{\alpha~\beta~\rho} \sum_{s\leq L} 2^{-Ln} 2^{s\alpha}\Big[2^{\rho s}+2^{\ell+L}\Big]^{\beta-m}~\leq~
\\\\ \ds
\C_{\alpha~\beta~\rho}\sum_{s\leq (\ell+L)/\rho} 2^{-Ln}2^{s\alpha}2^{(\ell+L)(\beta-m)}~+~\C_{\alpha~\beta~\rho}\sum_{ (\ell+L)/\rho<s\leq L} 2^{-Ln} 2^{s\alpha}2^{\rho s(\beta-m)}
\\\\ \ds
\leq~\C_{\alpha~\beta~\rho} ~2^{-Ln}2^{(\ell+L)\big[\alpha+\rho\beta-\rho m\big]/\rho}~+~\C_{\alpha~\beta~\rho}~2^{L\big[\alpha+\rho\beta-(n+\rho m)\big]}.
\end{array}
\eeq
From (\ref{Case4 Est1}) and (\ref{Case4 Est5}), we have
\bel{Case4 Est6}
\begin{array}{lr}\ds
\iint_{\Q_{0\ell}} \left|\I_{\alpha\beta}^\rho a(x,y)\right|^q dxdy~\leq~
\\\\ \ds
\C_{\alpha~\beta~\rho} \iint_{\Q_{0\ell}} 2^{-qLn}2^{q(\ell+L)\big[\alpha+\rho\beta-\rho m\big]/\rho}dxdy+\C_{\alpha~\beta~\rho} \iint_{\Q_{0\ell}} 2^{qL\big[\alpha+\rho\beta-(n+\rho m)\big]}dxdy
\\\\ \ds
\leq~\C_{\alpha~\beta~\rho} ~2^{-qLn}2^{q(\ell+L)\big[\alpha+\rho\beta-\rho m\big]/\rho} 2^{Ln} 2^{(\ell+L)m}~+~\C_{\alpha~\beta~\rho}~2^{qL\big[\alpha+\rho\beta-(n+\rho m)\big]}2^{Ln} 2^{(\ell+L)m}
\\\\ \ds
=~\C_{\alpha~\beta~\rho} ~2^{-L(q-1)n} 2^{(\ell+L)\big[q(\alpha+\rho\beta)-(q-1)\rho m\big]/\rho}~+~\C_{\alpha~\beta~\rho}~2^{L\big[q(\alpha+\rho\beta)-(q-1)n-q\rho m\big]} 2^{(\ell+L)m}.
\end{array}
\eeq
Note that ${\alpha+\rho\beta\over n+\rho m}=1-{1\over q}\Longrightarrow q(\alpha+\rho\beta)=(q-1)(n+\rho m)$ which further shows $q(\alpha+\rho \beta)>(q-1)\rho m$.
By using (\ref{Case4 Est6}) and summing over $\ell\ge0\colon \ell<(\rho-1)L$, we obtain
\bel{Case4 Est7}
\begin{array}{lr}\ds
\sum_{\ell\ge0\colon \ell<(\rho-1)L} \iint_{\Q_{0\ell}} \left|\I_{\alpha\beta}^\rho a(x,y)\right|^q dxdy~\leq~
\\\\ \ds
\C_{\alpha~\beta~\rho} ~\sum_{ \ell<(\rho-1)L}2^{-L(q-1)n} 2^{(\ell+L)\big[q(\alpha+\rho\beta)-(q-1)\rho m\big]/\rho}~+~\C_{\alpha~\beta~\rho}~\sum_{ \ell<(\rho-1)L} 2^{L\big[q(\alpha+\rho\beta)-(q-1)n-q\rho m\big]} 2^{(\ell+L)m}
\\\\ \ds
~=~\C_{\alpha~\beta~\rho~q}~ 2^{L\big[q(\alpha+\rho\beta)-(q-1)(n+\rho m)\big]}~=~\C_{\alpha~\beta~\rho~q}.
\end{array}
\eeq

 {\small School of Mathematical Sciences, Fudan University, Shanghai, 200433, China}\\
{\small email: 20114010001@fudan.edu.cn}

{\small  Department of Mathematics, Westlake University, Hangzhou, 310014, China}\\
{\small email: wangzipeng@westlake.edu.cn}
\end{document}